\documentclass{article}

\usepackage{amsmath} 
\usepackage{graphicx}
\usepackage{grffile} 
\usepackage{hyperref}
\usepackage[backend=biber]{biblatex} 
\usepackage[margin=1in]{geometry} 
\begin{document}

\title{Variational Solvers for Irreversible Evolutionary Systems}
\author{Andrés A {León Baldelli}\footnote{$\partial$'Alembert Institute, CNRS, Sorbonne Universités, Place Jussieu 75252 Paris Cedex 05, France}, %
Pierluigi Cesana\footnote{Institute of Mathematics for Industry, Kyushu University, 744 Motooka, Nishi-ku, Fukuoka, 819-0395, Japan.}}
\date{\today}

\maketitle

\section{Summary}\label{summary}

We study irreversible evolutionary processes with a general energetic
notion of stability. We dedicate this contribution to releasing three
nonlinear variational solvers as modular components (based on
FEniCSx/dolfinx) that address three mathematical optimisation problems.
They are general enough to apply, in principle, to evolutionary systems
with instabilities, jumps, and emergence of patterns which is
commonplace in diverse arenas spanning from quantum to continuum
mechanics, economy, social sciences, and ecology. Our motivation
proceeds from fracture mechanics, with the ultimate goal of deploying a
transparent numerical platform for scientific validation and prediction
of large scale natural fracture phenomena. Our solvers are used to
compute \emph{one} solution to a problem encoded in a system of two
inequalities: one (pointwise almost-everywhere) constraint of
irreversibility and one global energy statement. As part of our
commitment to open science, our solvers are released as free software.

\section{Statement of need}\label{statement-of-need}

Quasi-static evolution problems arising in fracture are strongly
nonlinear \cite{marigo:2023-la-mecanique},
\cite{bourdin:2008-the-variational}. They can admit multiple solutions,
or none \cite{leon-baldelli:2021-numerical}. This demands both a
functional theoretical framework and practical computational tools for
real case scenarios. Due to the lack of uniqueness of solutions, it is
fundamental to leverage the full variational structure of the problem
and investigate up to second order, to detect nucleation of stable modes
and transitions of unstable states. The stability of a multiscale system
along its nontrivial evolutionary paths in phase space is a key property
that is difficult to check: numerically, for real case scenarios with
several length scales involved, and analytically, in the
infinite-dimensional setting. The current literature in computational
fracture mechanics predominantly focuses on unilateral first-order
criteria, systematically neglecting the exploration of higher-order
information for critical points.

To fill this gap, our nonlinear solvers offer a flexible toolkit for
advanced stability analysis of systems which evolve with constraints.

\section{Three solvers}\label{three-solvers}

\texttt{HybridSolver} (1) \texttt{BifurcationSolver,} (2) and
\texttt{StabilitySolver} (3) implement the solution of three general
purpose variational problems:

\begin{enumerate}
\def\labelenumi{\arabic{enumi}.}
\item
  A constrained variational inequality; that is first order necessary
  conditions for unilateral equilibrium.
\item
  A singular variational eigen-problem in a vector space; that is a
  bifurcation problem indicating uniqueness (or lack thereof) of the
  evolution path.
\item
  A constrained eigen-inequality in a convex cone; originating from a
  second order eigenvalue problem indicating stabilty of the system (or
  lack thereof).
\end{enumerate}

These numerical tools can be used to study general evolutionary problems
formulated in terms of fully nonlinear functional operators in spaces of
high or infinite dimension. In this context, systems can have surprising
and complicated behaviours such as symmetry breaking bifurcations,
endogenous pattern formation, localisations, and separation of scales.
Our solvers can be extended or adapted to a variety of systems described
by an energetic principle (or unilateral stability law, see \cite{}
below).

We exploit the solvers to attack the following abstract problem which
encodes a selection principle: 
\[
P(0):\text{ Given } T >0, \text{ find an } \text{ irreversible-constrained evolution } y_t\]
\[y_t: t\in [0, T]\mapsto X_t  
\text{ such that}\]
\[\text{[Unilateral Stability]} \qquad E(y_t) \leq E(y_t + z), \quad \forall z \in V_0 \times K^+_0\qquad [1]
\]

Above, \(T\) defines a horizon of events. The system is represented by
its total energy \(E\) and \(X_t\) is the time-dependent space of
admissible states. A generic element of \(X_t\) contains a macroscopic
field that can be externally driven (or controlled, e.g.~via boundary
conditions) and an internal field (akin to an internal degree of order).
In the applications of fracture, the kinematic variable is a
vector-valued displacement \(u(x)\) and the degree of order
\(\alpha(x)\) controls the softening of the material. Irreversibility
applies to the internal variable, hence an irreversible-constrained
evolution is a mapping parametrised by \(t\) such that \(\alpha_t(x)\)
is non-decreasing with respect to \(t\). Remark that the test space for
the internal order parameter \(K^+_0\) only contains positive fields
owing to the irreversibility constraint. The main difficulties are to
correctly enforce unilateral constraints and to account for the changing
nature of the space of variations.

\subsection{Software}\label{software}

Our solvers are written in \texttt{Python} and are built on
\texttt{DOLFINx}, an expressive and performant parallel distributed
computing environment for solving partial differential equations using
the finite element method \cite{dolfinx2023preprint}. It enables us
wrapping high-level functional mathematical constructs with full
flexibility and control of the underlying linear algebra backend. We use
PETSc \cite{petsc-user-ref}, petsc4py \cite{dalcinpazklercosimo2011},
SLEPc.EPS \cite{hernandez:2005-slepc}, and dolfiny \cite{Habera:aa} for
parallel scalability.

Our solver's API receives an abstract energy functional, a user-friendly
description of the state of the system, its associated constraints, and
the solver's parameters. Solvers can be instantiated calling

\begin{verbatim}
solver = {Hybrid,Bifurcation,Stability}Solver(
      E,              # An energy (dolfinx.fem.form) 
      state,          # A dictionary of fields describing the system
      bcs,            # A list of boundary conditions
      [bounds],       # A list of bounds (upper and lower) for the state 
      parameters)     # A dictionary of numerical parameters
\end{verbatim}

where \texttt{\cite{}} are required for the \texttt{HybridSolver},
and used calling \texttt{solver.solve(\textless{}args\textgreater{})}
which triggers the solution of the corresponding variational problem.
Here, \texttt{\textless{}args\textgreater{}} depend on the solver (see
the documentation for details).

\texttt{HybridSolver} solves a (first order) constrained nonlinear
variational inequality, implementing a two-phase hybrid strategy which
is \emph{ad hoc} for energy models typical of applications in damage and
fracture mechanics. The first phase (iterative alternate minimisation)
is based on a de-facto \emph{industry standard}, conceived to exploit
the (partial, directional) convexity of the underlying mechanical models
\cite{bourdin:2000-numerical}. Once an approximate-solution enters the
attraction set around a critical point, the solver switches to perform a
fully nonlinear step solving a block-matrix problem via Newton's method.
This guarantees a precise estimation of the convergence of the
first-order nonlinear problem based on the norm of the (constrained)
residual.

\texttt{BifurcationSolver} is a variational eigenvalue solver which uses
SLEPc.EPS to explore the lower part of the spectrum of the Hessian of
the energy, automatically computed performing two directional
derivatives. Constraints are accounted for by projecting the full
Hessian onto the subspace of inactive constraints
\cite{jorge-nocedal:1999-numerical}. The relevance of this approach is
typical of systems with threshold laws. Thus, the \texttt{solve} method
returns a boolean value indicating whether the restricted Hessian is
positive definite. Internally, the solver stores the lower part of the
operators' spectrum as an array.

\texttt{StabilitySolver} solves a constrained variational eigenvalue
inequality in a convex cone, to check whether the (restricted) nonlinear
Hessian operator is positive therein. Starting from an initial guess
\(z_0^*\), it iteratively computes (eigenvalue, eigenvector) pairs
\((\lambda_k, z_k)\) converging to a limit \((\lambda^*, z^*)\) (as
\(k\to \infty\)), by implementing a simple projection and scaling
algorithm \cite{moreau:1962-decomposition},
\cite{pinto-da-costa:2010-cone-constrained}. The positivity of
\(\lambda^*\) (the smallest eigenvalue) allows to conclude on the
stability of the current state (or lack thereof), hence effectively
solving P(0). Notice that, if the current state is unstable
(\(\lambda^*<0\)), the minimal eigenmode indicates the direction of
energy decrease.

\begin{figure}
\centering
\includegraphics[width=.8\textwidth]{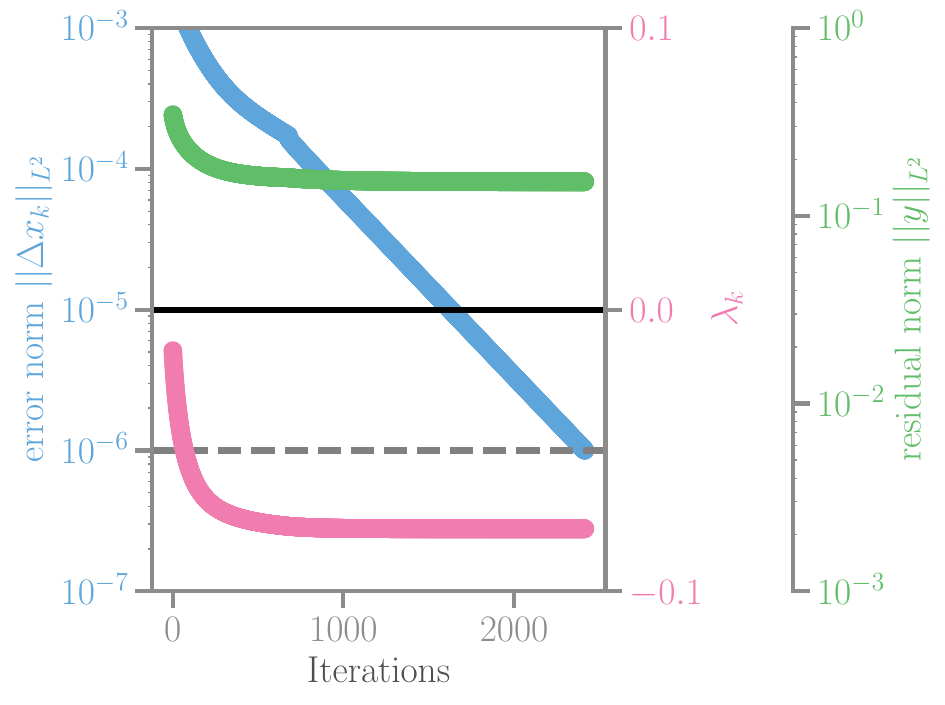}
\caption{Rate of convergence for \texttt{StabilitySolver} in 1d
(cf.~benchmark problem below). Targets are the eigenvalue
\(\lim_k \lambda_k =: \lambda^*\) (pink) and the associated eigen-vector
\(x^*\) (error curve in blue). Note that the residual vector (green) for
the cone problem need not be zero at a solution.\label{fig:convergence}}
\end{figure}

We dedicate a separate contribution to illustrate how the three solvers
are algorithmically combined to solve problem P(0) in the case of
fracture. \autoref{fig:convergence} illustrates the numerical
convergence properties of the \texttt{StabilitySolver} in the 1d
verification test.

\subsection{Verification}\label{verification}

\begin{figure}
\centering
\includegraphics[width=\textwidth]{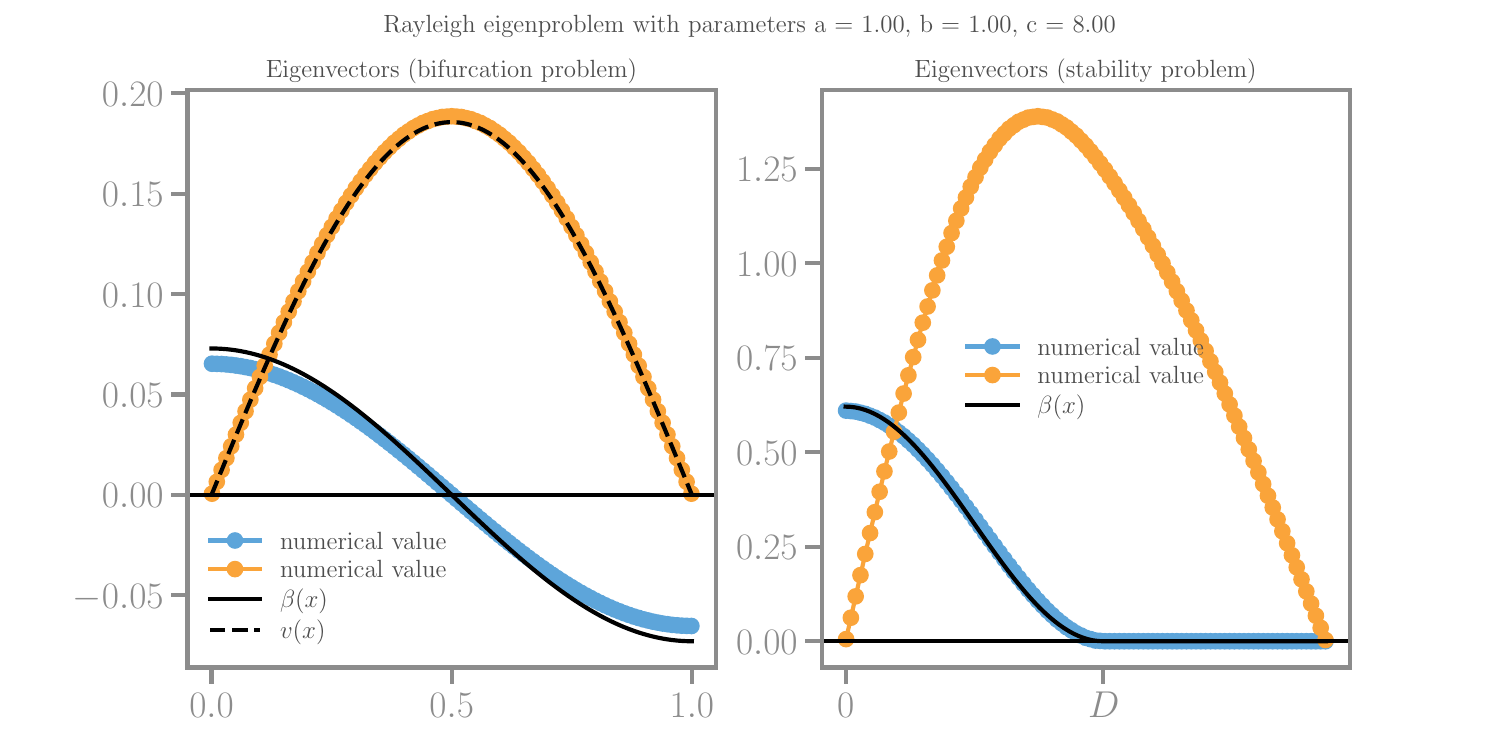}
\caption{Comparison between profiles of solutions \(\beta(x)\) in
\(X_0\) (left) vs.~\(\mathcal K^+_0\) (right). In the latter case, the
solution \(\beta(x)\) has support of size
\(D\in [0, 1]\).\label{fig:profile}}
\end{figure}

\begin{figure}
\centering
\includegraphics[width=\textwidth]{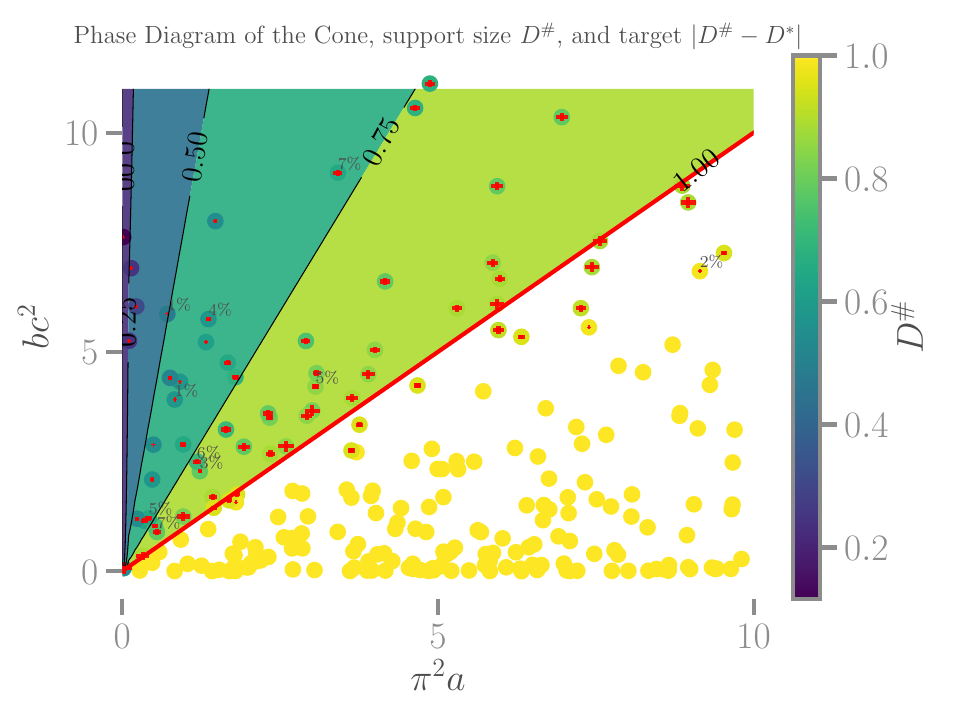}
\caption{The size of the support \(D\) for the minimiser in the cone.
Error bars indicate the absolute error. We observe a clear separation
between constant solutions with \(D=1\) vs.~nontrivial solutions. Where
error bars are not shown, the error is none.\label{fig:phase_diag_D}}
\end{figure}

\begin{figure}
\centering
\includegraphics[width=\textwidth]{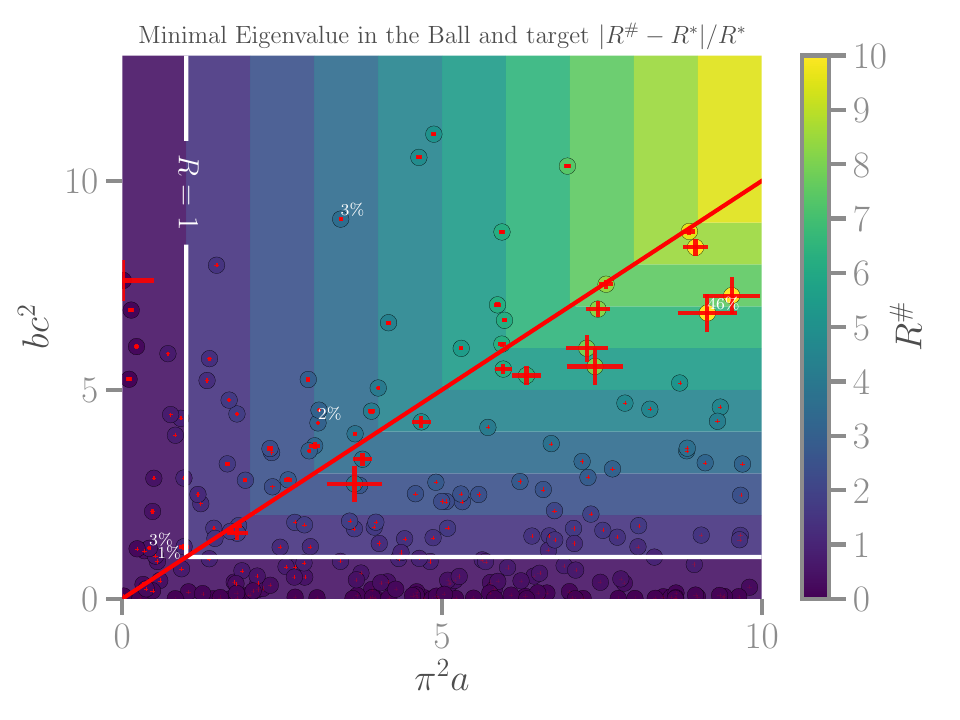}
\caption{Minimum value of \(\mathcal R\) in \(X_0\), numerical
computation vs.~closed form result. Notice the separation between
constant solutions and nontrivial solutions.\label{fig:phase_diag_ball}}
\end{figure}

\begin{figure}
\centering
\includegraphics[width=\textwidth]{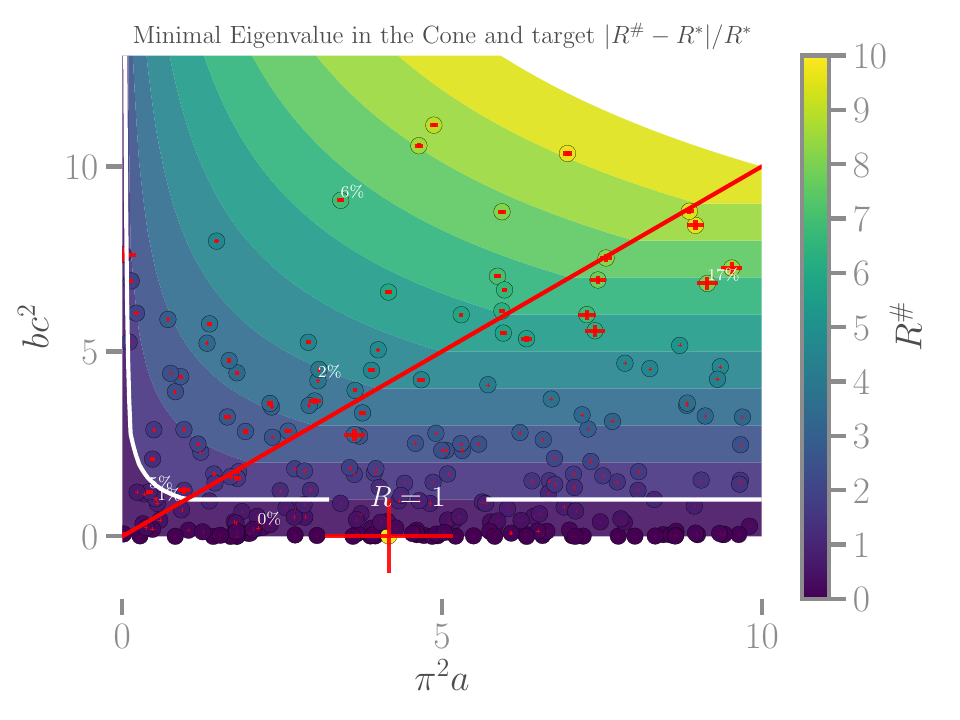}
\caption{Minimum value of \(\mathcal R\) in \(\mathcal K^+_0\),
numerical computation vs.~closed form results. The outlier
\((\pi^2a, bc^2)\sim (4, 0)\) represents a computation which did not
reach convergence. The mechanical interpretation is that only states
with \(R>1\) are energetically stable.\label{fig:phase_diag_cone}}
\end{figure}

We benchmark our solvers against a nontrivial 1d problem
(cf.~\texttt{test/test\_rayleigh.py} in the code repository), namely we
compute \[
\min_{X_0} \mathcal R(z) \quad \text{and} \quad \min_{\mathcal K^+_0} \mathcal R(z) \qquad\qquad [2],\]
using \texttt{BifurcationSolver} and \texttt{StabilitySolver}. The
quantity \(\mathcal R(z)\) is a Rayleigh ratio, often used in structural
mechanics as a dimensionless global quantity (an energetic ratio of
elastic and fracture energies) which provides insight into the stability
and critical loading conditions for a structure. For definiteness, using
the Sobolev spaces which are the natural setting for second order PDE
problems, we set \(X_0 = H^1_0(0, 1) \times H^1(0, 1)\) and
\(\mathcal K^+_0 = H^1_0(0, 1) \times \{\beta \in H^1(0, 1), \beta \geq 0\}\).
Let
\[\mathcal R(z):= \dfrac{\int_0^1 a(\beta'(x))^2dx+\int_0^1 b(v'(x) -c\beta(x))^2dx}{\int_0^1\beta(x)^2dx},\]
where \(a, b, c\) are real coefficients such that \(a>0, b>0, c\neq 0\).
The quantity above occurs in the stability analysis of a 1d damageable
bar, where \(b\) is related to the spring constant of the material while
\(a, c\) encapsulate material, loading, and model parameters, cf.~the
Appendix of \cite{pham:2011-the-issues}. \autoref{fig:profile} and
\autoref{fig:phase_diag_D} compare numerical results with the analytic
solution, \autoref{fig:phase_diag_ball} and
\autoref{fig:phase_diag_cone} show the relative error on the minimum in
the space of parameters.

\subsection{Acknowledgements}\label{acknowledgements}

A.L.B. acknowledges the students of MEC647 (Complex Crack Propagation in
Brittle Materials) of the
\texttt{Modélisation\ Multiphysique\ Multiéchelle\ des\ Matériaux\ et\ des\ Structures}
master program at ENSTA Paris Tech/École Polytechnique for their
contributions, motivation, and feedback; Yves Capdeboscq, Jean-Jacques
Marigo, Sebastien Neukirch, and Luc Nguyen, for constructive discussions
and one key insight that was crucial for this project. P.C. is a member
of the Gruppo Nazionale per l'Analisi Matematica, la Probabilità e le
loro Applicazioni (GNAMPA) of the Istituto Nazionale di Alta Matematica
(INdAM). P.C. holds an honorary appointment at La Trobe University and
is supported by JSPS Innovative Area Grant JP21H00102.

\subsection{References}\label{references}

\printbibliography

\end{document}